\documentclass[a4paper, english, 11pt]{article}

\usepackage{setspace}
\usepackage[T1]{fontenc}
\usepackage{hyperref}
\usepackage[utf8]{inputenc}
\usepackage{dcolumn}
\usepackage{color}
\usepackage[x11names]{xcolor}
\usepackage[thmmarks]{ntheorem}
\usepackage{ragged2e}
\usepackage{pdfpages}
\usepackage{amsmath,amssymb}
\usepackage[french, english]{babel}
\usepackage{multirow}
\usepackage{array}
\usepackage{booktabs}
\usepackage{wrapfig}
\usepackage{multicol}
\usepackage{algorithm}
\usepackage{lscape}
\usepackage{pifont}
\usepackage{graphicx}
\usepackage{cite}
\usepackage{calc}
\usepackage[all]{xy}
\usepackage{framed}
\usepackage{tikz}
\def\dar[#1#2]{\ar@<2pt>[#1]\ar@<-2pt>[#2]}

\usepackage{amsfonts,amssymb,amsmath}
\usepackage[top=3cm,bottom=3cm,left=3cm,right=3cm]{geometry}
\usepackage{fancyhdr}
\usepackage{array}

\newcommand{\kk}[0]{\textbf{k}}
\newcommand{\Mod}[0]{\text{Mod}}

\newcommand{\Pe}{\text{Pers}}

\newcommand{\R}[0]{\mathbb{R}}
\newcommand{\N}[0]{\mathbb{N}}

\newcommand{\Z}[0]{\mathbb{Z}}

\newcommand{\Obj}[0]{\text{Obj}}
\newcommand{\Ho}[0]{\text{H}}
\newcommand{\Hom}[0]{\text{Hom}}

\newcounter{nfigure} 
\makeatletter
 
\newskip\@bigflushglue \@bigflushglue = -100pt plus 1fil

\def\bigcentering{\let\\\@centercr\rightskip\@bigflushglue%
\leftskip\@bigflushglue
\parindent\z@\parfillskip\z@skip}

\theoremheaderfont{\scshape\bfseries}
\theorembodyfont{\normalfont}
\newtheorem{prop}{Proposition}[section] 

{ \theoremsymbol{\flushright{$\square$}} 
 \newtheorem*{pf}{Proof} }
 
  \newtheorem{rem}[prop]{Remarque} 
 
 \newtheorem{cor}[prop]{Corollary}
 \newtheorem{remark}[prop]{Remark}
 \newtheorem{lem}[prop]{Lemma}
 \newtheorem{defi}[prop]{Definition}
 \newtheorem{thm}[prop]{Theorem}
\author{Nicolas Berkouk}

\begin{document}
\selectlanguage{english}

\title{Stable resolutions of multi-parameter persistence modules}
\maketitle
\begin{abstract}
     The theory of persistence, which arises from topological data analysis, has been intensively studied in the one-parameter case both theoretically and in its applications. However, its extension to the multi-parameter case raises numerous difficulties. Indeed, it has been shown that there exists no complete discrete invariant for persistence modules with many parameters, such as so-called barcodes in the one-parameter case.
   
   To tackle this problem, some new algebraic invariants have been proposed to study multi-parameter persistence modules, adapting classical ideas from commutative algebra and algebraic geometry to this context. Nevertheless, the crucial question of the stability of these invariants has raised few attention so far, and many of the proposed invariants do not satisfy a naive form of stability. 
   
   In this paper, we equip the homotopy and the derived category of multi-parameter persistence modules with an appropriate interleaving distance. We prove that resolution functors are always isometric with respect to this distance, hence opening the door to performing homological algebra operations while "keeping track" of stability. This approach, we believe, can lead to the definition of new stable invariants for multi-parameter persistence, and to new computable lower bounds for the interleaving distance (which has been recently shown to be NP-hard to compute in \cite{BBK18}).
\end{abstract}

\tableofcontents

\section{Introduction}

Persistence theory appeared in the early 2000's as an attempt to make some constructions inspired by Morse theory actually computable in practice. For instance, in the context of studying the underlying topology of a data set. It has since been widely developed and applied. We refer the reader to \cite{Oudo15,Edel10} for extended expositions of the theory and of its applications.

The need for studying persistence modules obtained from functions valued in higher-dimensional partially ordered sets naturally arises from the context of data analysis, see for example \cite{Lesn15,LW15}.
However, as shown in \cite{Carl09}, the category $\Pe(\R^n)$ of functors $(\R^n,\leq) \to \Mod(\kk)$ seems to be too general at first sight for $n\geq 2$ to allow for some computer friendly analysis of its objects, as it contains a full sub-category isomorphic to the one of $\Z^n$-graded-$\kk[x_1,...,x_n]$-modules. 

Ideas from algebraic geometry \cite{Harr17} and combinatorial commutative algebra \cite{Ezr17} have since been developed in the context of persistence to tackle the study of the category of persistence modules with multiple parameters, from a computational point of view. Roughly speaking, those works propose some new invariants for summing up the algebraic structure of a given persistence module, arising from homological algebra methods. However informative about the algebraic structure of a given module they can be, the -crucial- question of the stability of these invariants is not studied yet. Indeed, data arising from the real world comes with noise, and one cannot avoid a theoretical study of the dependency of topological descriptors to small perturbations in view of practical applications. 

In this paper, we show that simple homological operations (as considering graded-Betti numbers of a persistence module) are not stable in a naive way with respect to the interleaving distance. To overcome this problem, we equip the homotopy (resp. derived) category of persistence modules of an interleaving-like distance and prove that the minimal free resolution functor (resp. localization functor) is distance-preserving with respect to the distances we introduce in section 5. 

Those distance comparison theorems open the door to perform homological algebra operations on free resolutions while \og keeping track of stability\fg. The future directions of research we shall undertake consist in studying precisely which meaningful stable invariants one can derive from those results. 

\subsection*{Content of the paper}

The paper is structured as follows : 

\begin{description}
\item[$\bullet$] Section 2 explains the motivating (counter)-example to naive stability of graded Betti numbers. We show that two persistence modules can be arbitrarily close in interleaving distance while having different graded Betti numbers. 

\item[$\bullet$] Section 3 aims at introducing precisely the category of persistence modules, and the classical notion of interleaving distance defined in the litterature.

\item[$\bullet$] Section 4 introduces some homological algebra background at the core of our constructions : the notions of homotopy and derived categories of an abelian category, together with projective resolution functors and localization functors. 

\item[$\bullet$] Section 5 introduces the new notions of interleavings in both case of the homotopy and the derived category of persistence modules. These induce two notions of interleaving distances on these two categories. We then prove our main theorems : the resolution and the localization functors are distance preserving with respect to these distances. As a corollary, we get a stability result for the Syzygy functor, from which are computed graded Betti numbers.

\item[$\bullet$] Section 6 goes back to the introducing example in the light of our results. We make explicit computations of homotopy interleavings that are not naive interleavings, hence understanding a form of stability for graded-betti numbers. We also discuss here future directions of work.
\end{description}

\subsection*{Acknowledgments}

The author would like to thank Steve Oudot and Magnus Botnan for many fruitful discussions and enlightnening comments on the subject of this article.

\section{Motivating Example : towards stable invariants for persistence with many parameters}

In this section, notions are not defined in a rigorous way, to allow the reader to build its intuition first. Persistence modules are algebraic structures that can encode the evolution of the topology along a filtration of a space. More precisely, they are functor from the category associated to the partially ordered set $(\R^n,\leq)$ (see section 3), to the category $\Mod(\kk)$ of vector spaces over a given field $\kk$. We can equip the category of persistence modules with the interleaving distance, and one can prove that this is the appropriate distance under which the topological variations (measured in $L_\infty$-norm) are transformed into algebraic variations for persistence modules. 

For the case where $n=1$, it has been proved in \cite{Craw12} that under weak assumptions on a persistence module $M : (\R,\leq) \to \Mod(\kk)$, it decomposes uniquely up to isomorphism as a direct sum of \emph{interval modules} : persistence modules that are constant on a given interval of $\R$ with value $\kk$ connected by internal maps being $\text{id}_\kk$. We record all the intervals that appear in the decomposition of a persistence module $M$ into a \emph{barcode}, and this invariant characterizes completely a persistence module over $\R$ up to isomorphism. Moreover, one can define the \emph{bottleneck distance} between barcodes, and prove that the interleaving distance between two modules is exactly the bottleneck distance between their barcodes (isometry theorem). Hence, the barcode of a persistence module is in particular a \emph{stable} invariant of persistence modules : a small perturbation with respect to the interleaving distance leads to a small perturbation of the barcode with respect to the bottleneck distance.  

As we shall see further, the algebraic structure of persistence modules for $n\geq 2$ is much more complicated to classify, and one shall not expect to get a nice (\emph{ie.} \og simple \fg) description of their isomorphism classes as it is the case when $n=1$. However, one can try to recover some partial informations about their structure by considering partial algebraic invariants. The first one proposed by Carlsson and Zormorodian in \cite{Carl09} is the $\emph{rank invariant}$, and records the ranks of the internal maps of a point-wise finite dimensional persistence module. This invariant generalizes in some sense the barcode of dimension $1$ case, as one can retrieve the barcode from the rank invariant of a one-parameter persistence module. Nevertheless, there are some very simple examples when $n\geq 2$ for which the rank invariant is not discriminative.

To strengthen this invariant, another proposition made computable in \cite{Chach16} are the \emph{graded Betti numbers} (see sections 2 and 5 of our paper). Their construction rely on the existence of minimal free resolutions for finitely generated persistence modules. Roughly speaking, taking a free resolution of a persistence module corresponds to choose a set of generators of the persistence module, then considering the relations between those generators and taking a set of generators of those relations, etc... Eventually this process stops (\emph{eg.} when the persistence module is finitely generated). Under some assumptions, one can choose at each step a \og minimal \fg ~set of generators and relations, leading to a minimal free resolution. We then define the $i$-th graded Betti number to be the function which associates to a point $x\in \R^n$ the number of generators appearing at $x$ in the $i$-th step of the free minimal resolution of a persistence module.  

We propose here to show that however informative about the algebraic structure of a persistence module, graded Betti numbers do not satisfy a naive form of stability. That is, we can find arbitrary close persistence modules with respect to the interleaving distance that have different Betti numbers.

Consider $\varepsilon \geq 0$, $M$ and $N_\varepsilon$ the persistence modules over $\R^2$ defined by : 

$$M = \kk_{\R_{\geq 0}\times \R_{\geq 0}}  ~~~\text{and} ~~~ N_\varepsilon = M \oplus \kk_{[0,\varepsilon]^2} $$

Where $\kk_A$ states for the persistence module with value $\kk$ on $A$ and internal maps between $s\leq t$ in $\R^2$ being $\text{id}_\kk$ if $s$ and $t$ are in $A$, the zero map otherwise.

\medskip 

\includegraphics[scale=1]{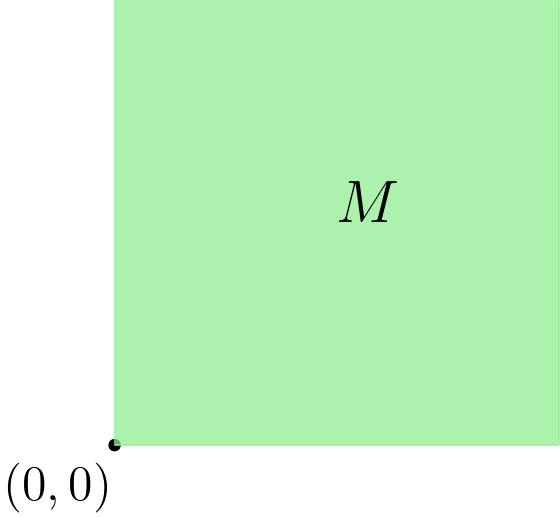} \hfill \includegraphics[scale=1]{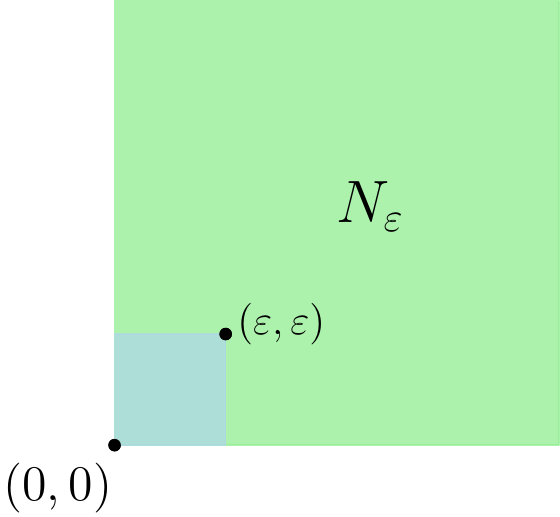}

Observe that the interleaving distance between $M$ and $N_\varepsilon$ is $\frac{\varepsilon}{2}$. Since $M$ is a free module (see definition in section 5), it is its own minimal free resolution. Hence, $\beta^0(M)(x) = 1$ for $x = (0,0)$ and $\beta^0(M)(x) = 0$ otherwise. 

Also, $N_\varepsilon$ has one more generator at $(0,0)$, thus $\beta^0(N_\varepsilon)(x) = 2$ for $x = (0,0)$ and $\beta^0(M)(x) = 0$ otherwise.

 Therefore, for any $\varepsilon \geq 0$ : 
 
 $$d_I(M,N_\varepsilon) =\frac{\varepsilon}{2} ~~~\text{and}~~~ \beta^0(M) \not =  \beta^0(N_\varepsilon) $$

This very simple example shows that Betti numbers are extremely sensitive to noise. In the following of the article, our aim will be to show how we can make stable the operation of \og taking resolution\fg, hence opening the door to obtaining new stable homological invariants for persistence with many parameters.

\section{Homological algebra in $\Pe(\R^n)$}

We fix a field $\kk$. Let $n >0$ an integer. We equip $\R^n$ with the product order, that is, for $(x_1,...,x_n)$ and $(y_1,...,y_n)$ in $\R^n$ : $$(x_1,...,x_n) \leq (y_1,...,y_n) \iff \forall ~1\leq i \leq n,~ x_i \leq_\R y_i  $$

Where $\leq_\R$ is the usual order on $\R$.

We look at $(\R^n,\leq)$ as a partially ordered set (poset) category, with objects elements of $\R^n$ and : $$\Hom_{(\R^n,\leq)}(x,y) = \begin{cases} \{*\} ~\text{if}~x\leq y \\
\emptyset ~\text{otherwise}

\end{cases} $$

\begin{defi}

A \textbf{persistence module} over $\R^n$ is a functor from $(\R^n,\leq)$ to the category of $\kk$-vector spaces $\Mod(\kk)$. We denote by $\Pe(\R^n)$ the category of persistence modules over $\R^n$ equipped with natural transformation of functors.
\end{defi}

\begin{remark} 
\begin{enumerate}
\item From its definition, it is clear that $\Pe(\R^n)$ is an abelian, complete and co-complete category. 

\item The main example of persistence modules arising from applications is the following. Consider $X$ a topological space and $f : X \to \R^n$ a map, then for $s\leq t$ in $\R^n$ and $i\in \Z_{\geq 0}$ we set :

\begin{itemize}

\item[$\bullet$] $M_i^f(s) = \Ho_i^\text{sing}(f^{-1}((-\infty, s]))$

\item[$\bullet$] $M_i^f(s\leq t) = \Ho_i^\text{sing}(f^{-1} ((-\infty, s])\subset f^{-1}((-\infty, t] ) )$
\end{itemize}

Where $\Ho_i^\text{sing}$ stands for the $i$-th singular homology with coefficient in $\kk$ and $(-\infty,t] = \{x\in \R^n \mid x \leq t\}$.

$M_i^f$ is the \textbf{$i$-th persistence module of $f$}.
\end{enumerate}

\end{remark}

To use $M_i^f$ as a descriptor of data coming from the real-world, we need to understand in which sense $M_i^f$ is sensitive to a perturbation of $f$. This is achieved by introducing the interleaving distance on $\Pe(\R^n)$. 

\begin{defi}
Let $M \in \Obj(\Pe(\R^n))$ and $\varepsilon \in \R$.
Define $M[\varepsilon]\in \Obj(\Pe(\R^n))$ the \textbf{$\varepsilon$-shift} of $M$ by, for $s\leq t $ in $\R^n$ : 
\begin{itemize}
\item[$\bullet$] $M[\varepsilon](s) = M \big( s+(\varepsilon,...,\varepsilon )\big) $

\item[$\bullet$] $M[\varepsilon](s\leq t) = M\big( s+(\varepsilon,...,\varepsilon) \leq t +(\varepsilon,...,\varepsilon) \big)$
\end{itemize}

\end{defi}

\begin{remark}
\begin{enumerate}

\item It is clear that the mapping $M \rightsquigarrow M[\varepsilon]$ induces an endofunctor of $\Pe(\R^n)$, called the \textbf{$\varepsilon$-shift} functor. It satisfies : $$\cdot[\varepsilon] \circ \cdot[-\varepsilon] = \cdot[-\varepsilon] \circ \cdot[\varepsilon] = \text{id}_{\Pe(\R^n)}$$

\item Most of the time, we will consider the case where $\varepsilon\geq 0$. In this situation, for $t\in \R^n$, define $s_\varepsilon^M(t) = M(t\leq t+(\varepsilon,...,\varepsilon))$. Then $s_\varepsilon^M$ is a natural transformation $M \Rightarrow M[\varepsilon]$ called the \textbf{$\varepsilon$-smoothing morphism}.
\end{enumerate}

\end{remark}

\begin{defi}
Let $M,N \in \Obj(\Pe(\R^n))$ and $\varepsilon \geq 0$. An \textbf{$\varepsilon$-interleaving} between $M$ and $N$ is the data of two morphisms $f : M \to N[\varepsilon]$ and $g : N \to M[\varepsilon]$ such that the following diagram commutes :

$$ \xymatrix{
M  \ar[rd]\ar@/^0.7cm/[rr]^{s_{2\varepsilon}^M} \ar[r]^{f} & N[\varepsilon]  \ar[rd] \ar[r]^{g[\varepsilon]} & M[2\varepsilon] \\
N  \ar[ur]\ar@/_0.7cm/[rr]_{s_{2\varepsilon}^N} \ar[r]^{g} & M[\varepsilon]  \ar[ur] \ar[r]^{f[\varepsilon]} & N[2\varepsilon]
   } $$ 
 
 Where the diagonal arrows are the smoothing morphisms.  
 If such a diagram exists, we say that $M$ and $N$ are \textbf{$\varepsilon$-interleaved} and write $M\sim_\varepsilon N$.
\end{defi}

\begin{remark}

Observe that $0$-interleaving corresponds to ismorphisms. Hence, one shall understand $\varepsilon$-interleaving as weaker form of isomorphisms. However, one must pay attention that \og being $\varepsilon$-interleaved\fg~ is not an equivalence relation since it lacks to be transitive. Indeed, if we have $M \sim_\varepsilon N \sim_\varepsilon O$ one can only deduce $M \sim_{2\varepsilon} O $.
\end{remark}

\begin{defi}
Let $M,N \in \Obj(\Pe(\R^n))$. Define their \textbf{interleaving distance} to be the possibly infinite number : $$d_I(M,N) = \inf\{\varepsilon \in \Z_{\geq 0} \mid M\sim_\varepsilon N\} $$
\end{defi}

\begin{prop}
The interleaving distance is an extended pseudo-distance on $\Pe(\R^n)$ that is, it satisfies for $M,N,O \in \Obj(\Pe(\R^n))$ : 
\begin{enumerate}
\item $d_I(M,N)\in \R_{\geq 0}\cup \{+\infty\}$
\item $d_I(M,N) = d_I(N,M)$
\item $d_I(M,N) \leq d_I(M,O) + d_I(O,N)$
\end{enumerate}
\end{prop}

The interleaving distance is the appropriate notion under which persistence modules are stable descriptors of real-world data in the following sense : 

\begin{thm}[Stability]
Let $X$ be a topological space, $f,g : X \to \R^n$ two maps and $i\in \Z_{\geq 0}$. Then : $$d_I(M_i^f,M_i^g) \leq \sup_{x\in X} |f(x) - g(x)|  $$
\end{thm}

Moreover, it has been shown in \cite{Lesn15} that it satisfies a form of universality property, justifying intrinsically the choice of this metric for $\Pe(\R^n)$.   

One way to understand the category $\Pe(\R^n)$, explained by Carlsson and Zormorodian in \cite{Carl09} in the specific case of persistence modules over $\N^n$, is to see persistence modules over $\R^n$ as $\R^n$-graded modules over the $\R^n$-graded algebra of generalized polynomials $\kk\{x_1,...,x_n\}$ (see definition below). This equivalence of category explains the complexity and impossibility to give a combinatorial classification of $\Pe(\R^n)$ when $n\geq 2$.

\begin{defi}
A \textbf{$\R^n$-graded $\kk$-algebra} is a $\kk$-algebra $A$ together with a decomposition of $\kk$-vector spaces $A = \bigoplus_{i\in\R^n}A_i$ such that $A_i \cdot A_j \subset A_{i+j}$ for all $i,j\in \R^n$.

Given $A$ a $\R^n$-graded $\kk$-algebra, a \textbf{$\R^n$-graded module} over $A$ is a module $M$ over the $\kk$-algebra $A$ together with a decomposition of $\kk$-vector spaces $M=\bigoplus_{i\in \R^n} M_i$ such that $A_i \cdot M_j \subset M_{i+j}$ for all $i,j\in \R^n$.

In both cases, $\R^n$-graded morphisms are usual morphisms that respect the decomposition. We denote by $A$-$\R^n$-grad-mod the abelian category of $\R^n$-graded modules over $A$.
\end{defi}

Let $\kk\{x_1,...,x_n\}$ be the algebra of generalized polynomials with coefficients in $\kk$ and real positive exponents, that is polynomials of the form $$P = \sum_i \alpha_i x_1^{d_i^1}... x_n^{d_i^n}, $$
with $\alpha_i \in \kk$ and $d_i^j\in \R_{\geq 0}^n$.

Observe that $\kk\{x_1,...,x_n\}$ is naturally a $\R^n$-graded $\R$-algebra with respect to the usual ring structure and decomposition $\kk\{x_1,...,x_n\}=\bigoplus_{d\in\R_{\geq 0}^n} \kk \cdot x^d$.

For $M$ a persistence module over $\R^n$, define $\alpha(M)$ to be the following $\kk\{x_1,...,x_n\}$-$\R^n$-graded module :

\begin{itemize}
\item[$\bullet$] $\R^n$-grading : $\alpha(M)=\bigoplus_{i\in\R^n}M(i)$
\item[$\bullet$] Action of $\kk\{x_1,...,x_n\}$ : for $d\in \Z_{\geq 0}^n$, define the action of $x^d$ component wise on $\alpha(M)$, that is for $i\in \R^n$, let $\cdot ~x^d : \alpha(M)_i = M(i) \to \alpha(M)_{i+d}=M(i+d) $  be the morphism $M(i\leq i+d)$
\end{itemize}

Conversely for $V = \bigoplus_{s\in\R^n} V_s$ a $\kk\{x_1,...,x_n\}$-$\R^n$-graded module, define $\beta(V)$ the persistence module over $\R^n$  by, for $s\leq t$ in $\R^n$ : 

\begin{itemize}
\item[$\bullet$] $\beta(V)(s) = V_s$
\item[$\bullet$] $\beta(V)(s\leq t)$ is the restriction of the action of $x^{t-s}$ to the component $V_s$
\end{itemize}

\begin{thm} \label{sec:moduleiso}
The mappings $\alpha$ and $\beta$ induce functors between $\Pe(\R^n)$ and $\kk\{x_1,...,x_n\}$-$\R^n$-grad-mod. These functors are additive exact isomorphisms of categories, inverse of each other.  
\end{thm}

\begin{rem}
\begin{enumerate}
     \item As an immediate corollary of theorem \ref{sec:moduleiso}, we deduce that the category $\Pe(\R^n)$ has enough projectives, that is, for any $M\in \Obj(\Pe(\R^n))$, there exists an epimorphism $P \twoheadrightarrow M$ with $P$ projective.
    \item With this theorem, we see that the study of persistence modules over $\R^n$ for $n\geq 2$ is a difficult problem, however, we can use standard tools of combinatorial homological algebra to perform computations.
   
\end{enumerate}

\end{rem}

\begin{defi}
Let $M\in \Obj(\Pe(\R^n))$. A \textbf{projective}  \textbf{resolution} of $M$ is a chain complex of objects of $\Pe(\R^n)$ : 

$$(F^\bullet, \partial_{F^\bullet}) :... \longrightarrow F^{-2} \stackrel{\partial_F^{-2}}{\longrightarrow} F^{-1} \stackrel{\partial_F^{-1}}{\longrightarrow} F^0 \stackrel{\partial_F^0}{\longrightarrow} 0 \longrightarrow ...$$

Such that : 

\begin{enumerate}

\item $F^j$ is projective for $j\leq 0$ and $F^j = 0$ for $j>0$
\item $\Ho^0(F^\bullet) := \dfrac{\mbox{Ker}(\partial^0)}{\mbox{Im}(\partial^{-1})}\simeq M$
\item $\Ho^j(F^\bullet):= \dfrac{\mbox{Ker}(\partial^j)}{\mbox{Im}(\partial^{j-1})} \simeq 0$ for $j \not = 0$
\end{enumerate}

\end{defi}

The existence of projective resolution for any $M\in \Obj(\Pe(\R^n))$ is a classical consequence from the fact that $\Pe(\R^n)$ have enough projectives.

\section{Homotopy and derived category}

In this section, we study the functorial properties of \og taking projective resolutions\fg. Let $A$ be an abelian category, $f : M_1\to M_2$ a morphism in $A$. Suppose that $M_1$ and $M_2$ admit projective resolutions $P_1^\bullet$ and $P_2^\bullet$, is there a chain map $\tilde{f}^\bullet : P_1^\bullet \to P_2^\bullet$ that lifts $f$ ? (ie. $\Ho^0(\tilde{f}^\bullet) = f$) ? Is it possible to make this lift functorial ?

To answer these questions, we will introduce the homotopy category of projective objects of $A$. A natural follow up then, is the construction of the derived category of $A$, which appears to be the right setting for homological computations in numerous areas of mathematics.

\subsection{Homotopy category and resolution functor}

For a general introduction to homotopy category, we refer the reader to \cite{Weib94,Oppe16}.
Let $A$ be an abelian category. We note $C(A)$ the category of chain complexes of $A$, and $C^-(A)$ its full sub-category with objects the complexes bounded from above, that is complexes $X^\bullet$ such that there exists $N\in \Z$ satisfying $X^n = 0$ for $n\geq N$.

\subsection{Chain homotopy and lifts}

\begin{defi}
Let $X$ and $Y$ be in $\Obj(C(A))$. An \textbf{homotopy} between $X$ and $Y$ is a chain map $\varphi : X\to Y$, such that there exists a collection $(h^i)_{i\in \Z}$ of morphisms in $A$, where for all $i \in \Z$, $h^i : X^i \to Y^{i-1}$ and  : $$ \varphi^i = \partial_Y^{i-1}\circ h^i + h^{i+1} \circ \partial_X^i $$

We note $\text{Hmt}(X,Y)$ the set of homotopies between $X$ and $Y$.  
\end{defi}

\begin{remark}
\begin{enumerate}
\item Observe that an homotopy induces the zero map at the level of cohomology objects : for all $i\in \Z$, $\Ho^i(\varphi) = 0$

\item $\text{Hmt}(X,Y)$ is a sub-group of $\Hom_{C(A)}(X,Y)$
\end{enumerate}

\end{remark}

\begin{defi}
With the same notations, let $f,g \in \Hom_{C(A)}(X,Y)$. We say that $f$ and $g$ are \textbf{homotopic} if $f-g\in \text{Hmt}(X,Y)$.

\end{defi}

\begin{remark}
\og Being homotopic \fg ~ is the equivalence relation associated to the subgroup $\text{Hmt}(X,Y)$, hence, it is compatible with the addition in $\Hom_{C(A)}(X,Y)$.
 \end{remark}

\begin{prop}
Let $\varphi \in \text{Hmt(X,Y)}$, $f\in \Hom_{C(A)}(W,X)$ and $g\in \Hom_{C(A)}(Y,Z)$. Then $g\circ \varphi\circ f \in \text{Hmt(W,Z)}$.
\end{prop}

As a consequence of this proposition, composition is well defined up to homotopy :

\begin{defi}
We define the \textbf{homotopy category} $K(A)$ of $A$ by :

\begin{itemize}
\item[$\bullet$] $\Obj(K(A)) = \Obj((C(A)))$
\item[$\bullet$] For $X,Y \in \Obj(K(A))$, $\Hom_{K(A)}(X,Y) =  \dfrac{\Hom_{C(A)}(X,Y)}{\text{Hmt(X,Y)}}$
\end{itemize}

\end{defi}

Until the end of the section, we assume that $A$ is an abelian category with enough projectives. Hence, every objects admits a projective resolution. 

\begin{prop}[Lift of projective resolutions]\label{prop:lift}
Let $X,Y \in \Obj(C(A))$, with projective resolutions respectively $P^\bullet$ and $Q^\bullet$. Let $f \in \Hom_{A}(X,Y)$, then : 
\begin{enumerate}
\item There exists a lift of $f$, that is a chain map ${f^\bullet}\in \Hom_{C(A)}(P^\bullet,Q^\bullet)$ such that  $\Ho^0(\tilde{f^\bullet}) = f$

\item For $g^\bullet \in \Hom_{C(A)}(P^\bullet,Q^\bullet) $ another lift, then $f^\bullet-g^\bullet \in \text{Hmt}(X,Y)$

\item For $P'^\bullet$ another projective resolution of $X$, there exists an isomorphism  $P^\bullet \simeq P'^\bullet$ in $K(A)$ 
\end{enumerate}
\end{prop}

Assume $A$ is a small category (that is, $\Obj(A)$ is a set), then by the axiom of choice, we can choose a projective resolution $P^\bullet(M)$ for each $M\in \Obj(A)$. From the above we get the following theorem : 

\begin{thm}
\begin{enumerate}
\item The association $M \rightsquigarrow P^\bullet(M)$  induces a well-defined functor $$P^\bullet : A \longrightarrow K(A)$$
\item $\Ho^0 \circ P^\bullet = \text{id}_A$
\item For any $j \not = 0$, $\Ho^j\circ P^\bullet = 0$
\end{enumerate}
\end{thm}

$P^\bullet$ is called a \textbf{resolution functor}. One shall remark that resolution functors are unique up to natural isomorphism, from proposition \ref{prop:lift}.

\subsection{Derived category and localization functor}

In this section, we give a very brief introduction to the derived category of an abelian category. Our  exposition follows the notes by Steffen Opperman \cite{Oppe16} and to which the reader can refer  for a detailed exposition of the construction. The idea behind the approach is to build a category in which we can make an object $M \in \Obj(A)$, seen as a complex concentrated in degree 0, isomorphic to its projective resolution $P^\bullet(M)$, in a precise sense.

\begin{defi}
Let $X,Y\in \Obj(C(A))$ and $f\in\Hom_{C(A)}(X,Y)$. Then $f$ is said to be a \textbf{quasi-isomorphism} if for all $j\in\Z$, $\Ho^j(f)$ is an isomorphism in $A$.

In the following, we will write qis to indicate quasi-isomorphisms.
\end{defi}

\begin{defi}
Let $X,Y\in \Obj(C(A))$. A \textbf{roof} from $X$ to $Y$ is the data of a chain complex $\tilde{X}\in \Obj(C(A))$, a quasi-isomorphism $q : \tilde{X} \to X $ and a chain morphism $f : \tilde{X} \to Y$, as summarized in the following diagram :

$$ \xymatrix{
   &  \tilde{X}  \ar[dl]_{q}^{qis} \ar[rd]^f &   \\
    X &  & Y
  } $$
  
For simplicity, we will write this roof $f\cdot q^{-1}$.
\end{defi}

\begin{lem}[Ore condition]
Let $X, \tilde{X}, Y \in \Obj(C(A))$, then given the solid part in the diagram below, it is possible to find the dashed part including $\tilde{Y}\in \Obj(C(A))$, such that the full diagram commutes in $C(A)$ :

$$ \xymatrix{
     \tilde{X}  \ar[d]_{q}^{qis} \ar[r]^f &Y \ar@{-->}[d]_{qis}^r   \\
    X \ar@{-->}[r] _g&  \tilde{Y} 
  } $$
  
 Dually, given the dashed part, we can find the solid part. 
\end{lem}

Two roofs $f_1\cdot q_1^{-1}$ and $f_2\cdot q_2^{-1}$ from $X$ to $Y$ are said to be \textbf{equivalent} if there exists a complex $H\in \Obj(C(A))$ and two quasi-isormophisms $h_i : H \to \tilde{X}_i$ ($i = 1,2$) such that the following commutes : 

$$ \xymatrix{
   &  \tilde{X_1}  \ar[dl]_{q_1} \ar[rd]^{f_1} &   \\
    X & H \ar[u]_{h_1}^{qis} \ar[d]^{h_2}_{qis} & Y\\
    &  \tilde{X_2}  \ar[ul]^{q_2} \ar[ur]_{f_2} &   \\
  } $$

\begin{prop}
Equivalence of roofs is an equivalence relation. Moreover, it is compatible with composition.

Where composition of roofs up to equivalence is given as follows : let $X,Y,Z \in \Obj(D(A))$, $f\cdot q^{-1}$ a roof from $X$ to $Y$, and $g\cdot r^{-1}$ from $Y$ to $Z$, define its composition to be the equivalence class of the following roof $\tilde{f}\cdot \tilde{r}^{-1}$ from $X$ to $Y$ given by the Ore condition lemma : 

$$ \xymatrix{
 & & \tilde{\tilde{X}} \ar@{-->}[ld]^{qis}_{\tilde{r}} \ar@{-->}[rd]^{\tilde{f}}& & \\
& \tilde{X} \ar[ld]^{qis}_q \ar[rd]^f & &\tilde{Y} \ar[ld]^{qis}_r \ar[rd]^g & \\
X & &Y & & Z
  } $$

\end{prop}

This ensures us that the following is well-defined: 

\begin{defi}
The \textbf{derived category} $D(A)$ of $A$ is given by : 

\begin{itemize}
\item[$\bullet$] $\Obj(D(A))= \Obj(C(A))$
\item[$\bullet$] For $X,Y\in \Obj(D(A))$, $\Hom_{D(A)}(X,Y) = \dfrac{\{\text{roofs from $X$ to $Y$}\}}{\text{equivalence of roofs}}$
\end{itemize}

\end{defi}

To each object $X\in \Obj(A)$ we can associate the corresponding chain complex concentrated in degree 0, and to each morphism $f : X \to Y \in \Hom_A(X,Y)$ we can associate the equivalence class of the roof $f \cdot \text{id}_X^{-1}$. This defines the fully faithful \textbf{localization functor} $\iota_A : A \to D(A)$.

We can do the exact same construction for $K(A)$, and define the \textbf{inclusion functor} $\iota_{K(A)} :  K(A) \to D(A)$ which is the identity on objects and sends homotopy class of morphisms $[f]\in\Hom_{K(A)}(X,Y)$ to the equivalence class of the roof $f \cdot \text{id}_X^{-1}$.

\begin{thm} \label{sec:thmdiag}
Assume $A$ has enough projectives and note $P^\bullet : A \to K(A)$ a resolution functor. Then the following diagram commutes up to natural isomorphisms :

$$ \xymatrix{
   &  A  \ar[dl]_{P^\bullet} \ar[rd]^{\iota_A} &   \\
    K(A)\ar[rr]_{\iota_{K(A)}}  &  & D(A)
  } $$

\end{thm}

We can actually refine the above. Observe that $P^\bullet$ actually defines a functor from $A$ to $K^-(\text{proj}-A)$ : the full sub-category of $K(A)$ with objects complexes $X$ such that there exists $N\in \Z$ satisfying $X^n=0$ for $n\geq N$, and $X^n$ is projective for $n < N$. Similarly, $\iota_{K(A)}$ restricts and co-restricts to a functor $\iota_{K^-(\text{proj}-A)} : K^-(\text{proj}-A) \to  D^-(A)$, with  $D^-(A)$ the full-subcategory of $D(A)$ with objects right bounded complexes.

\begin{thm}
The following diagram commutes up to natural isomorphism : 

$$ \xymatrix{
   &  A  \ar[dl]_{P^\bullet} \ar[rd]^{\iota_A} &   \\
    K^-(\text{proj}-A))\ar[rr]_{\iota_{K(\text{proj}-A)}}^\sim   &  & D^-(A)
  } $$
  
 Moreover, $\iota_{K(\text{proj}-A)}$ is an equivalence of categories.
\end{thm}

\section{Homotopy and derived interleavings}

In this section, we introduce the proper notions of interleaving at the levels of $K(\Pe(\R^n)))$ and $D(\Pe(\R^n))$ (hence inducing interleaving distances)  such that all arrows in the diagram of theorem \ref{sec:thmdiag} become isometries when taking $A$ to be $\Pe(\R^n)$.

\subsection{Homotopy interleaving distance}

\begin{defi}
Let $X \in \Obj(K(\Pe(\R^n)))$, and $\varepsilon \in \R_{\geq 0}$. Define the \textbf{$\varepsilon$-shift} of $X$, still written $X[\varepsilon ]$, to be the complex with $[\varepsilon]$ applied degree wise, that is : 

$$X[\varepsilon] = ~~~...\stackrel{\partial_X^{-2}[\varepsilon]}{\longrightarrow} X^{-1}[\varepsilon] \stackrel{\partial_X^{-1}[\varepsilon]}{\longrightarrow} X^{0}[\varepsilon] \stackrel{\partial_X^{0}[\varepsilon]}{\longrightarrow} X^{1}[\varepsilon] \stackrel{\partial_X^{1}[\varepsilon]}{\longrightarrow} ...$$

\end{defi}

\begin{remark}

\begin{enumerate}

\item Since $\cdot[\varepsilon]$ preserves homotopies, the association $X \rightsquigarrow  X[\varepsilon]$ induces a well defined endo-functor of $K(\Pe(\R^n)))$ 

\item Denoting $s^\varepsilon_X$ the collection of morphisms $(s^\varepsilon_{X^i})_{i\in \Z}$, it induces a chain morphism from $X \to X[\varepsilon]$, called the \textbf{$\varepsilon$-smoothing} morphism. We will write $[s^\varepsilon_X]$ for its homotopy equivalence class.

\end{enumerate}
\end{remark}
\begin{defi}[Homotopy interleaving]
Let $X,Y \in \Obj(K(\Pe(\R^n)))$. A \textbf{$\varepsilon$-interleaving} between $X$ and $Y$ is the data of two morphisms of $K(\Pe(\R^n))$, $f : X \to Y[\varepsilon]$ and $g : X \to Y[\varepsilon]$ such that the following diagram commutes in  $K(\Pe(\R^n)))$  :

$$ \xymatrix{
X  \ar[rd]\ar@/^0.7cm/[rr]^{[s_{2\varepsilon}^X]} \ar[r]^{f} & Y[\varepsilon]  \ar[rd] \ar[r]^{g[\varepsilon]} & X[2\varepsilon] \\
Y  \ar[ur]\ar@/_0.7cm/[rr]_{[s_{2\varepsilon}^Y]} \ar[r]^{g} & X[\varepsilon]  \ar[ur] \ar[r]^{f[\varepsilon]} & Y[2\varepsilon]
   } $$ 
   
 If such a diagram exists, we say that $X$ and $Y$ are \textbf{$\varepsilon$-interleaved} and write $X\sim_\varepsilon^K Y$.
\end{defi}

\begin{defi}
Let $X,Y \in \Obj(K(\Pe(\R^n)))$. Define their \textbf{homotopy interleaving distance} (interleaving distance when no confusion is possible) to be the possibly infinite number : 

$$d_I^K(X,Y) = \inf\{\varepsilon \in \Z_{\geq 0} \mid X\sim_\varepsilon^K Y\} $$
\end{defi}

\begin{prop}
The interleaving distance $d_I^K$ is an extended pseudo-distance on $K(\Pe(\R^n))$.
\end{prop}

\begin{pf}
The only tricky point here is to prove triangle inequality; which is an easy consequence from the fact that for $\varepsilon, \varepsilon '\in \Z_{\geq 0}$ and $X \in \Obj(K(\Pe(\R^n)))$: $$ \left[s_{\varepsilon '}^{X[\varepsilon]}\right]\circ \Big[s_\varepsilon^X\Big] = \Big[s_{\varepsilon + \varepsilon'}^X\Big] $$
\end{pf}

\subsection{Derived interleaving distance}
In this section, we equip the derived category of persistence modules over $\R^n$ with an interleaving distance, that will be proved to be appropriate later on with the derived comparison theorem.

We start by recalling a classical lemma of homological algebra stated as lemma 35.1 in \cite{Oppe16} :

\begin{lem}
Let $A,B$ two abelian categories and $F : A \to B$ an additive functor.

\begin{enumerate}
\item Applying $F$ degree wise gives rise to a well defined functor $F_K : K(A) \to K(B)$
\item There exists $F_D : D(A) \to D(B) $ such that the following diagram commutes up to natural isomorphisms if and only if $F$ preserves short exact sequences (ie. $F$ is \emph{exact}) :

$$\xymatrix{
    K(A) \ar[r]^{F_K} \ar[d]_{\iota_{K(A)}}  & K(B) \ar[d]^{\iota_{K(B)}} \\
    D(A) \ar[r]_{F_D} & D(B)
  }$$
  
  In this case, $F_D$ is naturally isomorphic to a functor acting on objects as $F$ applied degree wise. 

\end{enumerate}

\end{lem}

\begin{prop} \label{sec:shiftexact}
For any $\varepsilon \in \R$, the functor $\cdot [\varepsilon] : \Pe(\R^n) \to \Pe(\R^n)$ is exact.
\end{prop}

\begin{pf}
In $\Pe(\R^n)$, a sequence is exact if and only if it is exact point wise in $\Mod(\kk)$.
Let $M,N,O$ persistence modules fitting into the following short exact sequence : 

$$0 \longrightarrow M \stackrel{\varphi}{\longrightarrow} N \stackrel{\psi}{\longrightarrow} O \longrightarrow 0 $$

Then for any $s\in \R^n$, the following sequence is exact in $\Mod(\kk)$ : 

$$0 \longrightarrow M(s) \stackrel{\varphi_s}{\longrightarrow} N(s) \stackrel{\psi_s}{\longrightarrow} O(s) \longrightarrow 0 $$

And equivalently : 

$$0 \longrightarrow M(s+ \varepsilon) \stackrel{\varphi_{s+\varepsilon}}{\longrightarrow} N(s+\varepsilon) \stackrel{\psi_{s+\varepsilon}}{\longrightarrow} O(s+\varepsilon) \longrightarrow 0 $$

Which proves that the following sequence is exact : 

$$0 \longrightarrow M[\varepsilon] \stackrel{\varphi[\varepsilon]}{\longrightarrow} N[\varepsilon] \stackrel{\psi[\varepsilon]}{\longrightarrow} O[\varepsilon] \longrightarrow 0 $$\end{pf}

From the previous lemma, there is a well defined functor $[\varepsilon]_D : D(\Pe(\R^n)) \to D(\Pe(\R^n))$. For simplicity, we shall keep writing it $\cdot [\varepsilon]$ when no confusion is possible.

For $X \in \Obj(D(\Pe(\R^n)))$ and $\varepsilon \in \R_{\geq 0}$, we shall write $\{s_\varepsilon^X\}$ for $$\iota_{K(\Pe(\R^n))}([s_\varepsilon^X]) \in \Hom_{D(\Pe(\R^n))}(X,X[\varepsilon])$$

\begin{defi}[Derived interleaving]
Let $X,Y \in \Obj(D(\Pe(\R^n)))$. A \textbf{$\varepsilon$-interleaving} between $X$ and $Y$ is the data of two derived morphisms $f : X \to Y[\varepsilon]$ and $g : X \to Y[\varepsilon]$ such that the following diagram commutes in  $D(\Pe(\R^n)))$  :

$$ \xymatrix{
X  \ar[rd]\ar@/^0.7cm/[rr]^{\{s_{2\varepsilon}^X\}} \ar[r]^{f} & Y[\varepsilon]  \ar[rd] \ar[r]^{g[\varepsilon]} & X[2\varepsilon] \\
Y  \ar[ur]\ar@/_0.7cm/[rr]_{\{s_{2\varepsilon}^Y\}} \ar[r]^{g} & X[\varepsilon]  \ar[ur] \ar[r]^{f[\varepsilon]} & Y[2\varepsilon]
   } $$ 
   
 If such a diagram exists, we say that $X$ and $Y$ are \textbf{$\varepsilon$-interleaved} and write $X\sim_\varepsilon^D Y$.
\end{defi}

\begin{defi}
Let $X,Y \in \Obj(D(\Pe(\R^n)))$. Define their \textbf{derived interleaving distance} (or interleaving distance when no confusion is possible) to be the possibly infinite number : 

$$d_I^D(X,Y) = \inf\{\varepsilon \in \Z_{\geq 0} \mid X\sim_\varepsilon^D Y\} $$
\end{defi}

\begin{prop}
The interleaving distance $d_I^D$ is an extended pseudo-distance on $D(\Pe(\R^n))$.
\end{prop}

\subsection{Distance comparison theorems}

In this section, we compare the homotopy and derived interleaving distances we have introduced in the previous section with the usual interleaving distance. To do so, we look back at the commutative diagram of theorem \ref{sec:thmdiag}, with the additional structure of the appropriate notions of interleaving distances : 

$$ \xymatrix{
   &  (\Pe(\R^n),d_I) \ar[dl]_{P^\bullet} \ar[rd]^{\iota_{\Pe(\R^n)}} &   \\
    (K^-(\text{proj}-\Pe(\R^n))),d_I^K) \ar[rr]_{\iota_{K(\text{proj}-\Pe(\R^n))}}^\sim  &  & (D^-(\Pe(\R^n)),d_I^D)
  } $$

and we will prove that all functors appearing in the diagram are distance-preserving.

\subsubsection{Homotopy distance comparison theorem}
 Since $\Pe(\R^n)$ have enough projectives, we fix $P^\bullet$ to be a resolution functor.
 
\begin{prop}
Let $Q$ be a projective persistence module over $\R^n$, $M \in \Obj(\Pe(\R^n))$ and $\varepsilon \in \R$. Then the following holds : 
\begin{enumerate}
\item $Q[\varepsilon]$ is projective
\item There exists an isomorphism $P^\bullet(M[\varepsilon]) \simeq P^\bullet(M)[\varepsilon]$ in $K(\Pe(\R^n))$
\end{enumerate}
\end{prop}

\begin{pf}
\begin{enumerate}
\item We shall prove that the functor $\Hom_{\R^n}(Q[\varepsilon], -)$ is right exact (since it is always left exact). Consider an exact sequence in $\Pe(\R^n)$ : $$M \longrightarrow N \longrightarrow O \longrightarrow 0 $$

Then the sequence : $$M[-\varepsilon] \longrightarrow N[-\varepsilon] \longrightarrow O[-\varepsilon] \longrightarrow 0 $$

is still exact from proposition \ref{sec:shiftexact}. As $Q$ is projective, we get the following sequence in the category of abelian groups : 

$$\Hom_{\R^n}(Q,M[-\varepsilon]) \longrightarrow \Hom_{\R^n}(Q,N[-\varepsilon]) \longrightarrow \Hom_{\R^n}(Q,O[-\varepsilon]) \longrightarrow 0 $$

Since $\cdot[\varepsilon] \circ \cdot [-\varepsilon] = \text{id}_{\Pe(\R^n)}$, there are functorial isomorphisms for any $\mathcal{M} \in \Obj(\Pe(\R^n))$ :

$$\Hom_{\R^n}(Q,\mathcal{M}[-\varepsilon]) \simeq \Hom_{\R^n}(Q[\varepsilon],\mathcal{M}) $$

Leading to the desired exact sequence : 

$$\Hom_{\R^n}(Q[\varepsilon],M) \longrightarrow \Hom_{\R^n}(Q[\varepsilon],N) \longrightarrow \Hom_{\R^n}(Q[\varepsilon],O) \longrightarrow 0 $$

\item According to 1. and the fact that $\cdot [\varepsilon]$ is exact, $P^\bullet(M)[\varepsilon]$ is a projective resolution of $M[\varepsilon]$, hence it has to be isomorphic to $P^\bullet(M[\varepsilon])$ in $K(\Pe(\R^n))$.

\end{enumerate}
\end{pf}

\begin{prop} \label{sec:liftsmoothing}
Let $M\in \Obj(\Pe(\R^n))$ and $\varepsilon \geq 0$. Fix a chain morphism $\varphi : P^\bullet(M)[\varepsilon] \to P^\bullet(M[\varepsilon]) $ such that $[\varphi]$ is an isomorphism in $K(A)$. Then : $$P^\bullet(s_\varepsilon^M) =\left [\varphi \right ] \circ \left[s_\varepsilon^{P^\bullet(M)}\right]$$

Consequently : $$\Ho^0\left(\left[s_\varepsilon^{P^\bullet(M)}\right]\right) = s_\varepsilon^M $$
\end{prop}

\begin{pf}
Since  $\varphi \left ( P^\bullet(M[\varepsilon])\right ) = P^\bullet(M)[\varepsilon]$ it only remains to observe the commutativity of the following diagram in $C(\Pe(\R^n))$ :

$$\xymatrix{
      ... \ar[r]   & P^i(M) \ar[d]^{s^{P^i(M)}_\varepsilon} \ar[r] & ... \ar[r] & P^1(M) \ar[d]^{s^{P^1(M)}_\varepsilon} \ar[r] & P^0(M) \ar[r] \ar[d]^{s^{P^0(M)}_\varepsilon}  & M \ar[d]^{s^M_\varepsilon} \ar[r] & 0 \\
  ... \ar[r] & P^i(M)[\varepsilon] \ar[r] \ar[d]^{\varphi^i} & ... \ar[r] & P^1(M)[\varepsilon] \ar[r] \ar[d]^{\varphi^1} & P^0(M)[\varepsilon] \ar[r] \ar[d]^{\varphi^0} & M[\varepsilon] \ar[r] \ar[d]^{\text{id}_{M[\varepsilon]}}& 0 \\
   ... \ar[r] & P^i(M[\varepsilon]) \ar[r] & ... \ar[r] & P^1(M[\varepsilon]) \ar[r] & P^0(M[\varepsilon]) \ar[r] & M[\varepsilon] \ar[r] & 0
  } $$
  
  For simplicity, we have omitted the differential map of the chain complexes.
  
  Therefore, $\varphi \circ s_\varepsilon^{P^\bullet(M)}$ is one lift of $s_\varepsilon^M$, which by characterization of lift of morphism to projective resolutions proves that $P^\bullet(s_\varepsilon^M) = \left[\varphi \circ s_\varepsilon^{P^\bullet(M)}\right ] = [\varphi]\circ \left[s_\varepsilon^{P^\bullet(M)}\right]$.
\end{pf}

\begin{thm}[Homotopy comparison] \label{sec:hmtcompar}

Let $M,N\in \Obj(\Pe(\R^n))$ and $\varepsilon \in \R_{\geq}0$.
If $M$ and $N$ are $\varepsilon$-interleaved with respect to $f : M \to N[\varepsilon]$ and $g : N \to M[\varepsilon]$ in $\Pe(\R^n)$, then $P^\bullet(M)$ and $P^\bullet(N)$ are $\varepsilon$-interleaved with respect to $P^\bullet(f)$ and $P^\bullet(g)$ in $K(\Pe(\R^n))$.

Conversely, if $P^\bullet(M)$ and $P^\bullet(N)$ are $\varepsilon$-interleaved with respect to $\varphi :P^\bullet(M) \to  P^\bullet(N)[\varepsilon]$ and $\psi :P^\bullet(N) \to P^\bullet(M)[\varepsilon] $ in $K(\Pe(\R^n))$, then $M$ and $N$ are $\varepsilon$-interleaved with respect to $\Ho^0(\varphi)$ and $\Ho^0(\psi)$ in $\Pe(\R^n)$.

\end{thm}

\begin{pf}
Applying the functor $P^\bullet$ to the interleaving diagram gives : 

$$ \xymatrix{
P^\bullet(M)  \ar[rd]\ar@/^0.7cm/[rr]^{P^\bullet(s_{2\varepsilon}^M)} \ar[r]^{P^\bullet(f)} & P^\bullet(N)[\varepsilon]  \ar[rd] \ar[r]^{P^\bullet(g[\varepsilon])} & P^\bullet(M)[2\varepsilon] \\
P^\bullet(N)  \ar[ur]\ar@/_0.7cm/[rr]_{P^\bullet(s_{2\varepsilon}^N)} \ar[r]^{P^\bullet(g)} & P^\bullet(M)[\varepsilon]  \ar[ur] \ar[r]^{P^\bullet(f[\varepsilon])} & P^\bullet(N)[2\varepsilon]
   } $$ 
   
 Which is a $\varepsilon$-interleaving in $K(\Pe(\R^n))$ according to proposition \ref{sec:liftsmoothing}.
 
 The converse works similarly. 

\end{pf}

\begin{cor}
The functor $P^\bullet : (\Pe(\R^n),d_I) \to (K(\Pe(\R^n)),d_I^K) $ is distance preserving.
\end{cor}

\subsubsection{Derived distance comparison theorem}

\begin{thm}
\begin{enumerate}
\item The inclusion functor $\iota_{K(\Pe(\R^n))} : K(\Pe(\R^n)) \to D(\Pe(\R^n)) $ is interleaving preserving.

\item The localization functor $\iota_{\Pe(\R^n)} : \Pe(\R^n) \to D(\Pe(\R^n)) $ is interleaving preserving.
\end{enumerate}

\end{thm}

\begin{pf}
\begin{enumerate}
\item This is true by definition of derived interleavings
\item This is a consequence of the natural isomorphism $\iota_{K(\Pe(\R^n))}\circ P^\bullet \simeq \iota_{\Pe(\R^n)}$ and that isomorphisms define $0$-interleavings.
\end{enumerate}

\end{pf}

\begin{cor}
\begin{enumerate}
\item The inclusion functor $\iota_{K(\Pe(\R^n))} : (K(\Pe(\R^n)),d_I^K) \to (D(\Pe(\R^n)),d_I^D) $ is distance preserving. 
\item The localization functor $\iota_{\Pe(\R^n)} : (\Pe(\R^n),d_I) \to (D(\Pe(\R^n)),d_I^D) $ is distance preserving.
\end{enumerate}
\end{cor}

\subsection{Application to stable minimal free resolutions}

In this section, we use the theory of homotopy interleavings we have developped before to the special case of finitely generated persistence modules. In this setting, persistence modules have a canonical projective resolution that is particularly well behaved : the minimal free resolution. The graded-rank of each step of this resolution leads to define the graded Betti numbers of a finitely generated persistence module. By understanding the form of stability satisfied by \og taking minimal free resolution\fg, we hence explain why graded Betti numbers are not stable : they do not take homotopies into account. This phenomenom will be explicitely computed on an example in the next section.

Given $a \in \R^n$, we will denote by $F_a$ the persistence module over $\R^n$ defined by, for $s\leq t$ in $\R^n$ : 

$$ F_a(s) = \begin{cases} \kk ~~\text{if}~~a\leq s \\
0 ~~\text{else}

\end{cases} ~~~F_a(s\leq t) = \begin{cases} \text{id}_\kk ~~\text{if}~~a\leq s \\
0 ~~\text{else}

\end{cases} $$

\begin{defi} \label{sec:deffree}
A persistence module $M$ over $\R^n$ is \textbf{free} if there exists a function with finite support $\xi(M) : \R^n \to \Z_{\geq 0}$ such that : $$M \simeq \bigoplus_{i \in \R^n} F_i^{\oplus \xi(M)(i)}, $$ 

where $F_i^{\oplus \xi(M)(i)}$ corresponds to the direct sum of $\xi(M)(i)$ copies of $F_i$
\end{defi}

\begin{remark}
\begin{enumerate}

\item The name \og free \fg ~module corresponds to the notion of free modules introduced in \cite{EzrStu05} through $\alpha$ : $M$ is free in $\Pe(\R^n)$ if and only if $\alpha(M)$ is free in the sense of \cite{EzrStu05}. In particular, they are projective objects of the category $\Pe(\R^n)$.  
\item By the Krull-Schmidt theorem for finetely generated modules, the function $\xi(M)$ is unique for each $M$ since $F_a \simeq F_b $ if and only if $a=b$.

\end{enumerate}

\end{remark}

However we need some additional assumptions on $M$ to be assured of the existence of free resolutions. 

\begin{defi}
A persistence module over $\R^n$ is \textbf{finitely presented} if it is the quotient of two free persistence modules.
\end{defi}

It is clear that, up to rescaling, a finitely presented persistence module over $\R^n$ can be thought of as a discrete persistence module over $\Z^n$.
Also, adapting theorem \ref{sec:moduleiso} to this discrete case we get an additive isomorphism of categories between $\Pe(\Z^n)$ and the category $\kk[x_1,...,x_n]$-$\Z^n$-grad-mod, for which we have a similar definition of free modules.

\begin{defi}
Let $I = <x_1,...,x_n>$ be the maximal graded ideal of $\kk[x_1,...,x_n]$, $M\in \Obj(\Pe(\Z^n))$, and $(F^\bullet,\partial^\bullet)$ a projective (resp. free) resolution of $M$. Then $(F^\bullet,\partial^\bullet)$ is said to be \textbf{minimal} if for any $j \in \Z$, $\text{im} ~\alpha(\partial^j)\subset I \cdot \alpha(F_{j+1}) $. With $\alpha$ as in theorem \ref{sec:moduleiso}.

Now say that a free resolution of finitely presented $M\in \Obj(\Pe(\R^n))$ is minimal, if it is when discretized as a resolution of an object the associated object of $\Pe(\Z^n)$.
\end{defi}

The following classical theorem, can be found in chapter 8 of \cite{EzrStu05} in a slightly different form.

\begin{thm}[Hilbert's Syzygy]
Let $M\in \Obj(\Pe(\R^n))$ a finitely presented persistence module (ie. a quotient of a free persistence module). Then :

\begin{enumerate}
\item $M$ admits a minimal free resolution $(F^\bullet,\partial_F^\bullet)$
\item $F^i =0$ for $i<-n$, and $F^i$ is finitely generated for $-n\leq i \leq 0$
\item Any other minimal free resolution of $M$ is isomorphic to $(F^\bullet,\partial_F^\bullet)$ in $C(\Pe(\R^n))$.
\end{enumerate}

As a consequence, we can define, up to isomorphism, the minimal free resolution of $M$ and will denote it $(\pi^\bullet (M),\partial_{\pi^\bullet(M)})$.
\end{thm}

The previous theorem allows us to define the following : 

\begin{defi}
Let $M$ be a finitely generated persistence module over $\R^n$. Given $i\in\Z_{\geq 0}$, define its \textbf{$i$-th graded betti number} to be the function $\beta^i(M) : \R^n \to \Z_{\geq 0}$ such that :  $${\beta^i(M)= \xi (\pi^i(M))},$$

with $\xi$ defined in \ref{sec:deffree}. Observe that $\beta^i(M)$ is well defined since $\pi^i(M)$ is unique up to isomorphism.
\end{defi}

Let $\text{pers}(\R^n)$ be the full abelian subcategory of $\Pe(\R^n)$ with objects finitely generated modules. Then, by the Syzygy theorem, we can define the \textbf{Syzygy} functor $\pi^\bullet : \text{pers}(\R^n) \to K^-(\text{pers}(\R^n)) $ that associates to each finitely generated module its minimal free resolution. Then an immediate corollary of the homotopy distance comparison theorem \ref{sec:hmtcompar} is the following :

\begin{cor}
The functor $\pi^\bullet : (\text{pers}(\R^n),d_I) \to (K^-(\text{pers}(\R^n)),d_I^K) $ is distance preserving. 
\end{cor}

\begin{remark}
Minimal free resolutions of finitely presented persistence modules can be efficiently computed, following the work of W. Chacholski et al. \cite{Chach16}.
\end{remark}

\section{Computations and discussion}

\subsection{Back to our motivating example}
In this section, we develop the motivating example introduced in the beginning of the article in section 2. We get back to $M$ and $N_\varepsilon$ the persistence modules over $\R^2$. Using notations of \ref{sec:deffree}, the minimal free resolutions of $M$ and $N_\varepsilon$ are given, up to isomorphism, by the following : 

$$\xymatrix{
   \pi^\bullet(M)  \simeq 0 \ar[r] & 0 \ar[r] & 0 \ar[r] & F_{(0,0)} \ar[r] & 0 \\   
   \pi^\bullet(N_\varepsilon)  \simeq 0 \ar[r] & F_{(\varepsilon,\varepsilon)} \ar[r]^-{ \big( \begin{tiny} \begin{array}{c}
1\\
1 \\ \end{array} \end{tiny}\big)} & F_{(\varepsilon,0)} \oplus F_{(0,\varepsilon)} \ar[r]^-{ \big( \begin{tiny} \begin{array}{cc}
1 & 1\\
0 & 0 \\ \end{array} \end{tiny}\big)} & F_{(0,0)}^2 \ar[r] & 0 
}$$

Now, observe that for any $\eta > \frac{\varepsilon}{2}$,  $M$ and $N_\varepsilon$ are $\eta$-interleaved with respect to the following morphisms : 

$$\xymatrix{ M  \ar[rr]^{  \Big( \begin{footnotesize} \begin{array}{c}
0\\
s_\eta^{F_{(0,0)}} \\ \end{array} \end{footnotesize}\Big)  } & & N_\varepsilon [\eta] \\ 
N_\varepsilon \ar[rr]^{\big( 0 ~~ s_\eta^{F_{(0,0)}} \big) } & &  M [\eta]}$$

However it is clear that for every $i\in \{0,-1 -2\}$, $\pi^i(M)$ is not $\eta$-interleaved with $\pi^i(N_\varepsilon)$. Let us now construct a homotopy $\eta$-interleaving between $\pi^\bullet(M)$ and $\pi^\bullet(N_\varepsilon)$. 

Define the following complex : $$\xymatrix{
      C_\varepsilon = 0 \ar[r] & F_{(\varepsilon,\varepsilon)}  \ar[r]^-{ \big( \begin{tiny} \begin{array}{c}
1\\
1 \\ \end{array} \end{tiny}\big)} & F_{(\varepsilon,0)} \oplus F_{(0,\varepsilon)} \ar[r]^-{ ( \begin{tiny} \begin{array}{cc}
1 & 1\\
 \end{array} \end{tiny} )} & F_{(0,0)} \ar[r] & 0 
}$$

And observe that since $\pi^\bullet(N_\varepsilon) = \pi^\bullet(M) \oplus C_\varepsilon$, it is sufficient to prove that $C_\varepsilon$ is $\eta$-interleaved with $0$ in $K(\Pe(\R^n))$. This is equivalent to proving that $[s_{2\eta}^{C_\varepsilon}] = 0 $, that is, the chain map $s_{2\eta}^{C_\varepsilon} : C_\varepsilon \to C_\varepsilon[2\eta]$ is homotopic to $0$.

Now define : $$ \xymatrix{h^{-2} : C_\varepsilon^{-2} = F_{(\varepsilon,\varepsilon)} \ar[r]^-0 & 0 = C_\varepsilon^{-3}[2\eta] \\
h^{-1} : C_\varepsilon^{-1} =  F_{(\varepsilon,0)} \oplus F_{(0,\varepsilon)} \ar[r]^-{ \big( \begin{tiny} \begin{array}{c}
1\\
0 \\ \end{array} \end{tiny}\big)} &   F_{(\varepsilon - 2 \eta,\varepsilon - 2 \eta)} = C_\varepsilon^{-2}[2\eta] \\
h^{0} : C_\varepsilon^{0} =   F_{(0,0)} \ar[r]^-{ ( \begin{tiny} \begin{array}{cc}
0 & 1\\
 \end{array} \end{tiny} )} &  F_{(\varepsilon - 2\eta,-2\eta)} \oplus F_{(-2\eta,-2\eta)} = C_\varepsilon^{-1}[2\eta]
} $$

Note that $h^{-1}$ and $h^0$ are well-defined only when $\eta > \frac{\varepsilon}{2}$.
We have that $h=(h^i)$ defines an homotopy from $s_{2\eta}^{C_\varepsilon}$ to $0$. That is, we have  (with $h^i=0$ for $i$ different than $0,-1$) : $$s_{2\eta}^{C_\varepsilon} = d'\circ h + h \circ d $$

Where $d$ stands for the differential of the complex $C_\varepsilon$ and $d'$ the differential of $C_\varepsilon[2\eta]$. 

This example shows that, even in such a simple case, one cannot avoid taking into account homotopies in the problem of lifting interleavings to resolutions of persistence modules. Thus, to obtain homological invariants that are stable for persistence with many parameters, our work shows that a good algebraic framework is the homotopy category of persistence modules, equipped with the homotopy interleaving distance.

\subsection{Further directions of work}

\begin{itemize}

\item[$\bullet$] Our homotopy isometry theorem proves that, in order to interleave free resolutions, the important quantity to get interested in is the homotopy class of the smoothing morphism. Hence, there is no well defined notion of persistent image of a complex of persistence modules in the homotopy category. However, some derived techniques could be useful to overcome this difficulty, as studying some homological properties of the cone of this map, to give a criteria of closeness of two complexes of persistence modules. 

\item[$\bullet$] In \cite{BoLe16}, the authors prove that in the case of free persistence modules over $\R^2$, the interleaving distance can be computed exactly as a matching distance. This  result has since been extended in \cite{Bjer16}, where Bjerkevik proves an inequality bounding the bottleneck distance between two free persistence modules by a multilple of their interleaving distance. One could ask whether it is possible to define a bottleneck distance between minimal free resolutions of two persistence modules (that would allow to match free indecomposables across degrees), and to bound this distance by a multiple of the homotopy interleaving distance. This could lead to a computable lower bound to the interleaving distance (which has been shown to be NP-hard to compute for persistence modules with more than one parameter in  \cite{BBK18}), not relying on any kind of decomposition theorems.

\item[$\bullet$] In \cite{Kash18}, Kashiwara and Schapira have developed a sheaf theoretical framework for persistence, studied in detail in the case of dimension one in \cite{Berk18}. In many aspects, sheaf theory seems to be a convenient setting to deal with higher-dimensional persistence. However, its natural language is the one of the derived category. It is therefore necessary to enhance the derived category of persistence modules with the appropriate interleaving distance (which we have achieved in this paper), in order to compare constructions from both world.

\end{itemize}

\bibliographystyle{alpha}

\end{document}